\documentstyle[art10]{article}
\begin{document}

\newtheorem{lemma}{Lemma}
\newtheorem{theorem}{Theorem}
\newtheorem{definition}{Definition}
\newtheorem{corollary}{Corollary}
\newtheorem{proposition}{Proposition}

\def\R{{\bf R}}
\def\C{{\bf C}}
\def\Z{{\bf Z}}
\def\Q{{\bf Q}}
\def\A{{\cal A}}
\def\M{{\cal M}}
\def\vol{\mbox{vol}}

\title{\Large On nonformal simply connected symplectic manifolds}
\author{{\Large 
Ivan K. BABENKO 
\thanks{
    Department of Ma\-the\-ma\-tics and Mechanics, Moscow 
    State University,  
    Moscow, Russia, 
    e-mail: babenko@nw.math.msu.su, and
    D\'epartement de Math\'ematiques, Universit\'e 
    Montpellier II, Place Eug\`ene Bataillon,
    34095 Montpellier Cedex 5, France,  
    e-mail: babenko@math.univ-montp2.fr.}
\ and \ 
Iskander A. TAIMANOV
\thanks{
    Institute of Mathematics,
    630090 Novosibirsk, Russia,
    e-mail: taimanov@math.nsc.ru.
}
}}
\date{}
\maketitle

\begin{center}
{\bf 1. Introduction and main results}
\end{center}

A smooth manifold $M$ is called symplectic if it carries a nondegenerate
closed $2$-form $\omega$ which is called a symplectic form.
In this event a symplectic manifold means a pair $(M,\omega)$.
Since the skew-symmetric form $\omega$ is nondegenerate, 
$M$ is even-dimensional and moreover such a manifold always
carries an almost complex structure.

Topology of symplectic manifolds is being intensively developed now
(\cite{MS}) but still there are not many examples of compact symplectic
manifolds. The simplest ones are K\"ahler manifolds
and there are three constructions of new manifolds from
old ones. These are a symplectic fibration (\cite{Thurston}), a blow up
(\cite{McDuff}) and a fiber connect sum (\cite{Gompf}). The latter
two constructions were first outlined in \cite{Gromov}.

By results of Gromov (\cite{Nice}) and
Tischler (\cite{Tischler}), any
compact symplectic manifold is diffeomorphic to a  symplectic submanifold
of a complex projective space. The problem posed by Weinstein was to find
compact  symplectic manifolds not carrying a K\"ahler structure. The first
example of such a manifold was found by Thurston (\cite{Thurston}) and the
first simply connected example of such a manifold was constructed by McDuff
(\cite{McDuff}). Later Gompf constructed simply connected examples of the
minimal possible dimension which equals four (\cite{Gompf}).

An important property of K\"ahler manifolds is their 
formality established by
Deligne, Griffiths, Morgan, and Sullivan in \cite{DGMS}, which means that
the rational homotopy type is completely determined by 
the rational cohomology
ring. Formality of K\"ahler manifolds was used for distinguishing non
simply connected  symplectic manifolds with no K\"ahler structure in
\cite{BG,LO}.  For non simply connected spaces the notion of formality
is only defined for nilpotent spaces, whereas for
simply connected spaces it is always defined.
Meanwhile the problem of existence of nonformal simply connected
symplectic manifolds was open until now and moreover there was a conjecture
that such manifolds do not exist (the Lupton--Oprea conjecture, \cite{TO}).
We disprove it as follows

\vskip3mm

{\bf Main theorem.}
{\sl For any $N \geq 5$ there exist infinitely many  pairwise
non-homotopy-equivalent nonformal simply connected symplectic 
manifolds of
dimension $2N$.}

\vskip3mm

This theorem is proved as follows. Take a symplectic embedding of
the Kodaira--Thurston manifold $\widetilde{M}$
which is a symplectic four-dimensional nilmanifold (\cite{Thurston})
into $\C P^N$ for $N \geq 5$ and symplectically
blow up $\C P^N$ along $\widetilde{M}$ to obtain
a symplectic manifold $\widetilde{\C P}^N$.
It was McDuff who established that
$\widetilde{\C P}^5$ is a simply connected compact symplectic 
manifold with no
K\"ahler structure (\cite{McDuff}).
We prove that $\widetilde{\C P}^N$ is nonformal for $N \geq 5$
(see Theorem 2).

Moreover we construct an infinite series of
symplectic nilmanifolds $M(2m)$ of dimension $2m \geq 4$ and prove that
they are nonformal (Theorem 1) itself and that the blow ups of $\C P^N$ along
these manifolds embedded into $\C P^N$ are simply connected and
nonformal (Theorem 2).

Since the blow up at a point of a nonformal simply connected manifold is again
nonformal, successively applying blows up at points we obtain from one
nonformal  symplectic manifold infinitely many pairwise
non-homotopy-equivalent nonformal symplectic manifolds of the same
dimension.

In Section 7 we formulate and discuss some conjectures on nonformal
symplectic manifolds. We also would like to mention in 
the introduction that,
since all compact simply connected manifolds of dimension $\leq 6$ are
formal (\cite{NM}), the following problem stays open:

{\bf Problem.}
{\sl Are there nonformal simply connected symplectic manifolds
of dimension $8$ ?}

The results of the present paper in slightly weaker form
implying the existence of such manifolds of dimension $2N$ for any 
$N \geq 6$
were announced in \cite{BT}.

\vskip5mm

\begin{center}
{\bf 2. Some facts on symplectic manifolds}
\end{center}

Existence of an almost complex structure on symplectic manifolds is
demonstrated as follows.
Take a Riemannian metric $(\cdot, \cdot)$ on $M$ and define
an operator $A$ by the condition
$$
(Au,v) = \omega(u,v).
$$
The operator $A$ is skew-symmetric and hence the operator
$A^{\ast}A = - A^2$ is symmetric and positive. Take a positive
symmetric square root of it, $Q = \sqrt{-A^2} > 0$, and put
$J = AQ^{-1}$. It is clear that
$$
J^2 = -1,
$$
which means that $J$ defines an almost complex structure,
and
$$
\langle u, v \rangle = \omega(u,Jv)
$$
is a Hermitian metric on $M$, i.e., a Riemannian metric with respect to
which $J$ is skew-symmetric, which means that the almost complex structure
$J$ is compatible with the symplectic structure $\omega$.

This procedure enables us to introduce a smooth family of
compatible almost complex structures
on any smooth family of symplectic vector spaces, i.e., vector spaces with
symplectic forms.

Gromov proved that an open almost complex manifold $M$ always
carries a compatible symplectic structure (\cite{Nice}).

For compact manifolds existence of an almost complex structure does
not imply existence of a symplectic structure and the simplest
additional necessary condition is the
existence of a closed $2$-form $\omega$ such that its powers $\omega^j$
are cohomologically nontrivial for $j = 1,\dots,N$:
$[\omega]^j \neq 0$ in $H^{2j}(M)$.

A complex manifold $M$ is called a K\"ahler manifold if it carries
a Hermitian metric $h_{ij} dz^i d\bar{z}^j$ such that the form
$\omega = h_{ij} dz^i d\bar{z}^j$ is closed. This form is symplectic
and therefore any K\"ahler manifold
carries a natural  symplectic structure.

The simplest examples of K\"ahler manifolds are algebraic manifolds which
are complex submanifolds of the complex projective spaces. For such manifold
a K\"ahler structure is given by the metric induced 
from the Fubiny--Study metric
by the embedding. Denote by $(\C P^n, \omega_{FS})$
the complex projective space $\C P^n$ with a K\"ahler
form $\omega_{FS}$ induced by the Fubiny--Study metric. These symplectic
manifolds  serve as universal symplectic manifolds
in the following sense.

{\bf Proposition 1.}  (\cite{Nice, Tischler}) \
{\sl Let $(M,\omega)$ be a compact symplectic manifold of 
dimension $2n$ such that
the form $\omega$ is integer, i.e. $[\omega] \in H^2(M;\Z) \subset
H^2(M;\R)$. Then there exists an embedding
$$
f: M \rightarrow \C P^{2n+1}
$$
such that $f^{\ast}\omega_{FS} = \omega$.}

Rational forms, which are $\eta$ with $[\eta] \in H^2(M;\Q)$,
are dense in $H^2(M;\R)$ and hence
for any symplectic form $\omega$ on $M$ there
exists a small closed and nondegenerate perturbation $\omega +
\widetilde{\omega}$ of it such that the form 
$K(\omega + \widetilde{\omega})$ is
integer for some $K \in \R$ and 
therefore  $(M,K(\omega+\widetilde{\omega}))$
is a symplectic submanifold of $(\C P^{2n+1},\omega_{FS})$. Hence
symplectic submanifolds of the complex projective spaces
give us all topological types of compact symplectic manifolds.   

The first example of a compact symplectic 
manifold with no K\"ahler structure
is the Kodaira--Thurston manifold $\widetilde{M}$ which is defined 
as follows.

Denote by ${\cal H}$ the three-dimensional Heisenberg group formed
by the upper-triangular matrices
$$
\left(
\begin{array}{ccc}
1 & x & z \\
0 & 1 & y \\
0 & 0 & 1
\end{array}
\right)
$$
with $x, y, z \in \R$ and the multiplication operation.
The matrices with integer entries $x, y, z$ form a uniform lattice
${\cal H}_{\Z}$ in ${\cal H}$.

On the circle $S^1$, take a coordinate $u$ defined modulo $1$
and put
$$
\widetilde{M} = \left({\cal H} /{\cal H}_{\Z} \right) \times S^1
$$
with the form
\begin{equation}
\omega = dx \wedge du + dy \wedge dz.
\label{thurston}
\end{equation}
This is a four-dimensional symplectic nilmanifold whose one-dimensional real
cohomologies are generated by the forms $dx, dy$, and $du$,
and therefore the first Betti number equals three.
Since the odd-dimensional Betti numbers of K\"ahler manifolds are even,
$\widetilde{M}$ does not carry a K\"ahler structure but admits a complex
structure.

Later examples of four-dimensional symplectic nilmanifolds with
no complex structure were found (\cite{FGG}, see also Section 6).

\vskip5mm

\begin{center}
{\bf 3. On the blow up construction and simply 
connected symplectic manifolds
with no K\"ahler structure}
\end{center}

The blow up construction is defined for any pair $Y \subset X$ 
of smooth manifolds where the
structure group of the normal bundle to $Y$ reduces to $U(k)$ with
$2k  = \dim X - \dim Y$.

Let $(X,\omega)$ be a compact symplectic manifold of dimension $2N$ and
let $Y$ be a symplectic submanifold, of $X$, of dimension $2(N-k)$.
For any $p \in Y$ consider the subspace $E_p \subset T_pX$ which is
the orthogonal complement to $T_p Y$ with respect to $\omega$.
These subspaces form a bundle $E = E_Y$, over $Y$,
isomorphic to the normal bundle to $Y$. Since $Y$ is a symplectic
submanifold, we have

1) the restriction of $\omega$ to $E_p Y$ is nondegenerate for any
$p \in Y$;

2) the normal bundle to $Y$ in $X$ is naturally identified with $E$.

\noindent
By the method exposed in Section 2 construct a fiberwise almost complex
structure on $E$ compatible with the restrictions of $\omega$ to the 
fibres. This implies that the structure group of $E$ reduces to
$U(k) = SO(2k) \cap Sp(k)$.

Now let $Y$ be a submanifold of $X$ and the structure group
of the normal bundle $E$ to $Y$ be $U(k)$.
Identify the fibres of $E$ with $\C^k$ and
consider another bundle
$$
\widetilde{E} \rightarrow Y
$$
with fibres isomorphic to
the canonical line bundle over $\C P^{k-1}$.
This canonical line bundle is
$$
L \rightarrow \C P^{k-1}
$$
with
$$
L = \{(z,l) \in \C^k \times \C P^{k-1} | z \in l\}.
$$
The condition $z \in l$ is written as
$$
z_i l_j = z_j l_i \ \ \mbox{for $i,j = 1,\dots,k$}
$$
where $(l_1,\dots,l_k)$ are the homogeneous coordinates on $\C P^{k-1}$.
This bundle is associated with $E$ which means that the action of the
structure group on $\widetilde{E}$ is as follows
$$
A \cdot (z,l) = (Az, Al)
$$
where $A \cdot z = Az$ is the corresponding action of the structure group
of $E$.

The fibres of $E$ are endowed with a Hermitian metric and we denote by
$E_r$ and $\widetilde{E}_r$ the submanifolds of $E$ and $\widetilde{E}$
defined by the condition on $z$: $|z| \leq r$.

The fibre spaces
$\left(\widetilde{E}_1 \setminus \widetilde{E}_0 \right) \rightarrow Y$
and $\left( E_1 \setminus E_0 \right) \rightarrow Y$
are canonically isomorphic and their fibres are diffeomorphic to
a punctured disc $\{z \in \C^k | 0 < |z| \leq 1\}$.
The fibre space $E_0 \rightarrow Y$ is just 
the zero section of $E \rightarrow Y$, the fibre space
$$
\widetilde{Y} = \widetilde{E}_0 \rightarrow Y
$$
is called the projectivization of $E \rightarrow Y$, and its fibres are
diffeomorphic to $\C P^{k-1}$.

Now we are able to construct the blow up of $X$ along $Y$. For that take
a closed tubular neighborhood $V$ of $Y$ in $X$ and naturally identify it
with the fibre space $E_1$.
Now take the manifold
$$
\widetilde{X} = \overline{(X \setminus V)} \cup_{\partial V}
 \widetilde{E}_1 
$$
with $\widetilde{E}_1$ glued to the boundary of $\overline{(X \setminus V)}$
by the natural isomorphism $\partial E_1 = \partial \widetilde{E}_1$.

The manifold $\widetilde{X}$ is called the blow up of $X$ along $Y$ and
roughly speaking is obtained by replacing a closed neighborhood $V$ of
$Y$, which is a disc bundle over $Y$ by a manifold with boundary
$\widetilde{V}$ which is a disc bundle over $\widetilde{Y}$.

There is a natural projection
$$
\pi:\widetilde{X} \rightarrow X
$$
which is a diffeomorphism outside $Y$ and whose restriction to
$\pi^{-1}(Y)$ is the fibering map $\widetilde{Y} \rightarrow Y$.

Notice that $\partial V = \partial E_1 = \partial 
\widetilde{E}_1$ fibers
over $Y$ with the fibre $S^{2k-1}$ and $\widetilde{E}_1$ fibers
over $Y$ with the fibre $\overline{\C P^k \setminus D^{2k}}$ where
the overline stands for closure. The fibre of the latter bundle
is homotopy equivalent to $\C P^{k-1}$ and the embedding $$
\widetilde{i}: \partial V
\rightarrow \widetilde{E}_1 = \widetilde{V}
$$
extends to a commutative diagram
of the fibre bundles
\begin{equation}
\begin{array}{ccccccc}
\partial V = \partial \widetilde{E}_1 & &&
\longrightarrow & && \widetilde{E}_1 \\
& S^{2k-1}&\searrow & & \swarrow & \sim\C P^{k-1}& \\
& & & Y & & &
\end{array}
\label{triangle}
\end{equation}
and the horizontal mapping $\widetilde{i}$ preserves the fiberings.

The Leray--Hirsch theorem implies

{\bf Proposition 2.}
{\sl 1) The cohomology ring of the projectivization $\pi:
\widetilde{Y} \rightarrow Y$
of a vector bundle $E \rightarrow Y$ is isomorphic to
$$
H^{\ast}(\widetilde{Y}) = H^{\ast}(Y)[a]/\langle a^k + c_1 a^{k-1} 
+ \dots +
c_{k-1}a + c_k \rangle
$$
where $ c_1,\dots, c_k$ are the Chern classes of the vector bundle
$E \rightarrow Y$. The induced homomorphism $\pi^{\ast}:H^{\ast}(Y)
\rightarrow H^{\ast}(\widetilde{Y})$ is a monomorphism.

2) The cohomology class $a \in H^2(\widetilde{Y})$
may be chosen such
that $\widetilde{i}^{\ast}(a) = 0$ in $H^2(\partial V)$.}

The proof of the first statement of this proposition is exposed,
for instance, in \cite{Husemoller}. It may be explained also in terms of
spectral sequences of fiberings from (\ref{triangle}) and from the
functoriality of these sequences it follows that $a$, a generator of
$E^{0,2}_2 = H^0(Y;H^2(\C P^{k-1}))$, is mapped by $\widetilde{i}^{\ast}$ 
into $H^0(Y;H^2(S^{2k-1})) = 0$. This finishes the proof of the
proposition.

{\bf Proposition 3.} (\cite{McDuff}) \
{\sl If $k \geq 2$, then $\pi_1(\widetilde{X}) = \pi_1(X)$.}

The proof of this proposition is as follows. By the Van Kampen theorem,
$$
\pi_1(\widetilde{X}) = \pi_1(\overline{X \setminus V})
\ast_{\pi_1(\partial V)}
\pi_1(\widetilde{E}_1).
$$
Notice that $V$ retracts onto $Y$. Since $\dim X - \dim Y \geq 4$,
any $2$-disc immersed into $X \setminus Y$ is deformed into a disc with no
intersections with $Y$ and this implies that
 $\pi_1(\overline{X \setminus V}) = \pi_1(X)$.

Since both fibres in (\ref{triangle})
are simply connected, the projections induce isomorphisms
of the fundamental groups. The diagram is commutative and therefore
the embedding $\partial \widetilde{E}_1 \rightarrow \widetilde{E}_1$
also induces an isomorphism of the fundamental groups which together
with the Van Kampen proves the proposition.

The simplest example is the blow up along a point $p = Y$ 
in a $2n$-dimensional
manifold $X$ which results in adding the complex projective space:
$$
\widetilde{X} = X \# \overline{\C P^n}
$$
where the overline stands for the opposite orientation.

In  \cite{Gromov} Gromov outlined the blow up for symplectic manifolds
and the detailed exposition of it was done by McDuff (\cite{McDuff}):

{\bf Proposition 4.} (\cite{Gromov,McDuff})\
{\sl If $Y$ is a compact symplectic submanifold of a symplectic manifold
$(X,\omega)$ then the blow up $\widetilde{X}$ of $X$ along $Y$ 
carries a symplectic 
form $\widetilde{\omega}$
which equals $\pi^{\ast}\omega$ outside a neighborhood of $\pi^{-1}(Y)$.}

Before explaining the McDuff example let us 
recall two fundamental properties of
compact K\"ahler manifolds $M$:

1) for any $k$ the cohomology group $H^k(M;\C)$ has a decomposition
$$
H^k(M;\C) = \oplus_{p+q=k} H^{p,q}(M;\C)
$$
where the elements of $H^{p,q}(M;\C)$ are realized by forms of the type
$\alpha dz^{i_1} \wedge \dots \wedge dz^{i_p} \wedge d\bar{z}^{j_1} \wedge
\dots \wedge d\bar{z}^{j_q}$ and the complex conjugation establishes
isomorphisms $H^{p,q}(M;\C) \rightarrow H^{q,p}(M;\C)$
(the Hodge decomposition). This implies,
in particular, that if $k$ is odd then $\dim H^k(M;\C)$ is even;

2) if the complex dimension of $M$ is $n$ then the mapping
$$
H^{n-k}(M;\C) \stackrel{\times [\omega]^k}\longrightarrow H^{n+k}(M;\C)
$$
is an isomorphism for any $k=0,\dots,n$ (the Hard Lefschetz condition).

Now we are able to describe

{\bf The McDuff manifold.}

The form (\ref{thurston}) is integer and therefore,
by the Gromov--Tischler theorem (see Proposition 1), there exists
a symplectic embedding of the Kodaira--Thurston manifold
$\widetilde{M}$ into $\C P^5$. Define $\widetilde{\C P}^5$ as
the blow up of $\C P^5$
along $\widetilde{M}$.

{\bf Proposition 5.} (\cite{McDuff}) \
{\sl 1) The dimension of $H^3(\widetilde{\C P}^5;\C)$ equals $3$. 
Therefore $H^{\ast}(\widetilde{\C P}^5;\C)$ has no 
Hodge decomposition and
$\widetilde{\C P}^5$ has no K\"ahler
structure.

2) The symplectic manifold $\widetilde{\C P}^5$
does not satisfy the Hard Lefschetz condition.}

\vskip5mm

\begin{center}
{\bf 4. Minimal models and formality}
\end{center}

{\bf A) Differential graded algebras and their minimal models.}

Recall that a differential graded algebra 
is a graded algebra
$$
\A = \oplus_{k \geq 0} \A^k
$$
with a differential $d: \A \rightarrow \A$ of degree $1$
, i.e., $d(\A^k)) \subset \A^{k+1}$, such that

1) $x \wedge y = (-1)^{kl} y \wedge x $ for $x \in \A^k, y \in \A^l$;

2) $d(x \wedge y) = dx \wedge y + (-1)^k x \wedge d y$ for $x \in \A^k$
(the Leibnitz rule);

3) $d^2 = 0$.

This implies that the cohomology ring $H^{\ast}(\A)$ of $\A$ is defined
and this graded ring, supplied with the zero differential $d =0$, 
is also 
a differential graded algebra.

In the sequel we consider algebras over  the ground field $F$
which may be $\Q$  (the rational numbers), $\R$ (the real numbers), or 
$\C$
(the complex numbers) and in addition assume that $\dim_F \A^k$ is finite
for any $k$.

A homomorphism $f$ of two differential graded algebras $(\A,d_{\A})$ and
$({\cal B},d_{\cal B})$ is a homomorphism $f:\A \rightarrow {\cal B}$
of algebras such that
$$
f(\A^k) \subset {\cal B}^k \ \ \mbox{and} \ \ f(da) = df(a),
$$
i.e., it preserves the grading and the differential.
This implies that any such homomorphism induces a homomorphism
$$
f^{\ast}: H^{\ast}(\A) \rightarrow H^{\ast}({\cal B})
$$
of the cohomology rings. 

An algebra is called connected if $H^0(\A) = F$ where $F$ is the 
ground field and
it is called one-connected if in addition $H^1(\A) = 0$.

Let $x_1,x_2,\dots$ be a sequence of weighted variables with
$\deg x_i \geq 1$ for any $i$. Denote by $\Lambda(x_1,\dots)$ the
free graded commutative algebra generated by $x_1,\dots$.

A differential graded algebra $\M$ is called minimal if

1) $\M = \Lambda(x_1,\dots)$ for some family of free generators and
there are only finitely many generators of any fixed degree;

2) $dx_i \in \Lambda(x_1,\dots,x_{i-1})$.

Notice that for a one-connected algebra condition 2) may be  replaced by

2') $d(\M) \subset \M^+ \wedge \M^+$ where $\M^+ = \oplus_{k>0} \M^k$.

We say that $\M(\A)$ is a minimal model for $\A$ if

1) $\M(\A)$ is a minimal algebra;

2) there is a homomorphism $h:\M \rightarrow \A$ inducing an isomorphism of
the cohomology rings.

The fundamental theorem of Sullivan says

{\bf Proposition 6.} (\cite{Sullivan}) \
{\sl Every one-connected differential graded algebra has a minimal model
unique up to isomorphism.}

Now, for any simply connected polyhedron $X$ such that the ranks of 
its homotopy groups are finite, construct the 
minimal algebra $\M_X$ over $\Q$
as follows.   

Take the Postnikov tower of $X$ which
is a family of spaces $X_k$ and mappings $f_k:X \rightarrow X_k$
such that

1) $\pi_j(X_k) = 0$ for $j > k$;

2) $f_k$ induces isomorphisms $f_{k\ast}: \pi_j(X) 
\rightarrow \pi_j(X_k)$
for $j \leq k$;

3) there are commutative up to homotopy diagrams
$$
\begin{array}{ccc}
X & \stackrel{f_{k+1}}{\rightarrow} & X_{k+1} \\
 & f_k \searrow & \downarrow \\
 & & X_k
\end{array}
$$
where $p_{k+1}: X_{k+1} \rightarrow X_k$ is a fibration with fibre
$K(\pi_{k+1}(X),k+1)$.

By the Cartan--Serre theorem,

1) $H^{\ast}(K(\Z,n);\Q)$ is a free graded commutative algebra
$\Lambda(x)_n$, over $\Q$,
generated by $x$ with $\deg x = n$ which means that 
$$
\Lambda(x)_n = \Q[x] \ \mbox{for $n$ even and} \ \Lambda(x)_n =
\Q[x]/x^2 \ \mbox{for $n$ odd};
$$

2) $H^j(K(G,n);\Q) = 0$ for any finite commutative group $G$ 
and $j \geq 1$.

\noindent
This remains true after replacing $\Q$ by $\R$ or $\C$. 

Using the Cartan--Serre theorem together with the Postnikov tower, construct
$\M_X$ inductively.

For $X = K(\pi,n)$ with $n \geq 2$ put $M_X = H^{\ast}(K(\pi,n);\Q)$.

Assume that for $X_k$ the algebra $\M_k = \M_{X_k}$ is defined and
$H^{\ast}(\M_k;\Q) = H^{\ast}(X_k;\Q)$.
Consider the cohomology spectral sequence for the fibration $X_{k+1}
\rightarrow X_k$ and take the transgression homomorphism
$$
d_{k+1}: E^{k+1,0} = H^{k+1}(K(\pi_{k+1}(X),k+1);\Q) \rightarrow
E^{0,k+2} = H^{k+2}(X_k;\Q).
$$
Let $y_1,\dots,y_i$ be linear generators 
of $H^{k+1}(K(\pi_{k+1}(X),k+1);\Q)$
and for any $y_j$ take an element $w_j$ in $\M_k$ representing $d_{k+1}y_j$.
Now put $dy_j = w_j$ and define
$$
\M_{k+1} = \M_k \otimes_d
\Lambda(\mbox{Hom}\,(\pi_{k+1}(X),\Q))_{k+1}
$$
where
$$
\Lambda(\mbox{Hom}\,(\pi_{k+1}(X),\Q))_{k+1}
= H^{k+1}(K(\pi_{k+1}(X),k+1);\Q),
$$
$$
d:\mbox{Hom}\,(\pi_{k+1}(X),\Q) \rightarrow \M_k^{k+1},
$$
and the differential on $\M_{k+1}$ is determined by its restrictions
to $\M_k$ and $\Lambda$. Taking other representatives of $dy_j$ we
obtain isomorphic algebras.

Now consider the series of such extensions
$$
\dots \subset \M_k \subset \M_{k+1} \subset \dots \subset \M_{\infty},
$$   
take the limit algebra $\M_{\infty} = \cup_{k \geq 0}\M_k$
and put $\M_X = \M_{\infty}$.

It follows from the construction that there is a natural isomorphism
\begin{equation}
\mbox{Hom}\,(\pi_{\ast}(X),\Q) = \M_X / \M_X \wedge \M_X.
\label{spheric}
\end{equation}
Other important properties were also established by Sullivan who had
invented minimal models.

{\bf Proposition 7.} (\cite{Sullivan}) \
{\sl 1) For a simply connected compact polyhedron $X$
there is a differential algebra ${\cal E}(X)$ of $\Q$-polynomial forms
on $X$ and $\M_X$ is the minimal model for ${\cal E}(X)$;

2) The algebra $\M_M \otimes \R$ is the minimal model for the
algebra ${\cal E}^{\infty}(M)$ of smooth differential forms on a
compact manifold $M$.}

For these reasons, $\M_X$ is called the minimal model for $X$,
up to isomorphism it describes
the rational homotopy type of $X$.

An important property of minimal models which is based on their construction
via the Postnikov towers and is also partially justified by
(\ref{spheric}) and Proposition 7 is the following.

{\bf Proposition 8.}
{\sl Every simplicial mapping of polyhedra or smooth mapping of manifolds
$f:X \rightarrow Y$ induces a homomorphism of
the minimal models
$$
\hat{f}:\M_Y \rightarrow \M_X
$$
such that the induced homomorphisms $f^{\ast}:H^{\ast}(Y;\Q) \rightarrow
H^{\ast}(X;\Q)$ and $\hat{f}^{\ast}:H^{\ast}(\M_Y;\Q) \rightarrow
H^{\ast}(\M_X;\Q)$ coincide as
 $$f^{\ast}h^{\ast}_Y = h^{\ast}_X\hat{f}^{\ast},$$
where $h_X$ and $h_Y$ are maps of minimal models in ${\cal E}(X)$ and
${\cal E}(Y)$ correspondingly.} 

For a detailed exposition of minimal models see \cite{DGMS,GM} where
the minimal model is also defined for compact nilpotent polyhedra.
The nilpotence of a polyhedron means that its fundamental group is nilpotent
and its actions on higher homotopy groups are also nilpotent.

{\bf B) Minimal models for nilmanifolds.}

Recall that a nilmanifold $N = G/\Gamma$ is the compact quotient space
where $G$ is a simply connected nilpotent group and $\Gamma$ a uniform
lattice in $G$.

We formulate only a few facts about the minimal models for
nilmanifolds. Given a nilmanifold $X = G/\Gamma$,
there is
a uniquely defined (up to isomorphism) graded differential
algebra $\M= \M_X$ over $\Q$ satisfying, in particular,
the following conditions

1) $\M$ is freely generated by elements $x_1,\dots,x_k$ of degree 1;

2) $\M$ is a minimal algebra;

3) $H^{\ast}(\M_X) = H^{\ast}(X;\Q)$.

This algebra is called the minimal model for a nilmanifold $X$ and
has a nice and clear algebraic origin.

Let ${\cal G}$ be the Lie algebra of a nilpotent group $G$ and
${\cal G}^{\ast}$ be the dual algebra to ${\cal G}$. The Lie brackets
define a mapping
$$
{\cal G} \times {\cal G} \stackrel{[\cdot,\cdot]}{\rightarrow} {\cal G},
$$
and it is easily checked that the dual mapping
$$
d: {\cal G}^{\ast} \rightarrow {\cal G}^{\ast} \times
{\cal G}^{\ast}
$$
is a differential and in this event the equality
$d^2 = 0$ is equivalent to the Jacobi identity.
The nilpotence of ${\cal G}$ implies minimality of the algebra
$\Lambda({\cal G}^{\ast},d)$ which will be the minimal model for $G/\Gamma$.

The differential of this minimal model is written in terms of the structure
constants as follows. Let $\{e^1,\dots,e^k\}$ be a basis for ${\cal G}$ and
$\{\omega_1,\dots,\omega_k\}$ be the dual basis for ${\cal G}^{\ast}$.
Let in this basis the Lie brackets be given by
$$
[e^i,e^j] = \sum_k c^{ij}_k e^k.
$$
Then the differential $d$ is as follows
\begin{equation}
d \omega_k = \sum_{i,j} c^{ij}_k \omega_i \wedge \omega_j.
\label{diff}
\end{equation}

Elements of the Lie algebra are naturally associated
with left-invariant vector fields on $G$ and elements of ${\cal G}^{\ast}$
are represented in this event by left-invariant $1$-forms on $G$.

Assume that $G$ possesses uniform lattices. By the results of Mal'tsev
(\cite{Maltsev}) this happens if, and only if,
the structure constants
$c^{ij}_k$ are rational. Let $\Gamma \subset G$ be a uniform lattice
and $X = G/\Gamma$ be the corresponding nilmanifold. By the Nomizu theorem
(\cite{Nomizu}), there is a natural isomorphism
$$
H^{\ast}(X;F) = H^{\ast}(\widehat{\Omega}(1),F)
$$
where $\widehat{\Omega}(1)$ is a subcomplex generated as 
a graded differential
algebra by left-invariant $1$-forms on $G$ and the isomorphism is induced
by an embedding of the complexes of forms. For $F = \Q$ it 
needs to consider
forms with rational periods  and in this events the rationality of the
structure constants guarantees the ``rationality'' of the differential
(\ref{diff}).

We have

{\bf Proposition 9.}
{\sl Given an $n$-dimensional nilmanifold  $X = G/\Gamma$, 
its minimal model
$\M_X$ is correctly defined and is as follows:

1) $\M_X = (\Lambda(x_1,\dots,x_n);d)$ with $\deg x_i = 1$
for $i = 1,\dots,n$;

2) $ d x_k = \sum_{i,j} c^{ij}_k x_i \wedge x_j$.}

{\sl Example.}
The minimal model of the Kodaira--Thurston manifold $\widetilde{M}$ 
is freely generated by elements $\eta_1,\dots,\eta_4$
of degree 1 such that
\begin{equation}
d\eta_1 = d\eta_2 = d\eta_4 = 0, \ d\eta_3 = \eta_1 \wedge \eta_2.
\label{thurston2}
\end{equation}
These elements are realized via Proposition 7 by the following
left invariant forms on $\widetilde{M}$:
$$
\eta_1 = dx,  \ \eta_2 = dy, \ \eta_3 = xdy + dz, \ \eta_4 = du.
$$

By the Gompf theorem (\cite{Gompf}), any finitely presented group is the
fundamental group of a four-dimensional symplectic compact manifold and
therefore the minimal model is not well defined for all symplectic
manifolds.     

{\bf C) Formality of differential algebras and spaces.}

A homomorphism of differential graded algebras
$$
({\cal A},d_{\cal A}) \to ({\cal B},d_{\cal B})
$$
is called quasiisomorphism if it induces an isomorphism of cohomologies.

A minimal algebra ${\cal M}$ is called formal if there is a quasiisomorphism
$$
({\cal M},d) \to (H^{\ast}({\cal M}),0).
$$
In particular, this implies that $({\cal M},d)$ is the minimal model
for its cohomology ring $(H^{\ast}({\cal M}),0)$ with 
a zero differential.

A differential graded algebra ${\cal A}$ is called formal if its
minimal model is formal.

A sufficient condition for formality that is effective for applications
follows from Proposition 6:
if there is a quasiisomorphism 
from a one-connected differential graded algebra $({\cal A},d)$ to 
its cohomology ring $(H^{\ast}({\cal A}),0)$ or if it exists in the other 
direction then the algebra ${\cal A}$ is formal.

A polyhedron or a smooth manifold $X$
is called formal if its minimal model $\M_X$ is formal.
For such a space it is said that its rational homotopy type is a
formal consequence of its cohomology.

Among the examples of formal spaces there are Riemannian locally symmetric
spaces (\cite{GM}), classifying spaces (\cite{Sullivan}), compact 
K\"ahler
manifolds (\cite{DGMS}), and simply connected compact manifolds 
of dimension $\leq 6$ (\cite{NM}).

The simplest example of a nonformal manifold is the three-dimensional
non simply connected but nilpotent space ${\cal H}/{\cal H}_{\Z}$.

{\bf Proposition 10.}
{\sl The nilmanifold ${\cal H}/{\cal H}_Z$ is nonformal.}

By Proposition 9,
the minimal model $\M$ for ${\cal H}/{\cal H}_Z$ is as follows:
$\M = \Lambda(x_1,x_2,$
$x_3), \deg x_i = 1$, and $dx_1 = 
dx_2 = 0, dx_3 =
x_1 \wedge x_2$.
Assuming that this manifold is formal, there exists
a homomorphism $\psi:\M \rightarrow H^{\ast}(\M)$ inducing the identity of
cohomologies. Notice that $\psi(x_1) \cup \psi(x_2) \in
\psi(x_1) \cup H^1(\M) = 0$ which implies that
$\psi(x_1 \wedge x_3) = \psi(x_1) \cup \psi(x_3) = 0$.
The element $x_1 \wedge x_3$
is closed but not exact in $\M$ and hence
$\psi(x_1 \wedge x_3) \neq 0$.
We arrive at a contradiction which proves the proposition.

By (\ref{thurston2}), the same arguments show that the Kodaira--Thurston
manifold $\widetilde{M}$ is also nonformal.

An analysis of this situation leads to the following criterion.
The minimal model for a simply connected space $X$
is isomorphic as a graded commutative algebra
to
$$
\M = \otimes_{k \geq 0} \Lambda(V_k)_k
$$
where $V_k = \mbox{Hom}(\pi_k(X),k),\Q)$.
In each $V_k$ take the subspace $C_k$ of closed elements.

{\bf Proposition 11.} (\cite{DGMS})\
{\sl $\M$ is formal if, and only if,  in each $V_k$ there is 
a complement $N_k$
to $C_k$, $V_k = C_k \oplus N_k$, such that any closed 
form in the ideal $I_N$
generated by the elements of $N_k$, $I_N = I(\oplus N_k)$, is exact.}

{\bf D) Massey products.}

We define only the triple Massey product.

Let $a \in \M^p, b \in \M^q$, and $c \in \M^r$ represent 
nontrivial cohomology classes
such that $[a] \cup [b] = 0$ and $[b] \cup [c] =0$.
Therefore there exist $g \in \M^{p+q-1}$ and $h \in \M^{q+r-1}$
such that $a \wedge b = dg$ and $b \wedge c = dh$.

Define a cycle
$$
k = g \wedge c + (-1)^{p-1} a \wedge h.
$$
Its cohomology class is defined modulo 
$([a]H^{q+r-1}(\M) + [c]H^{p+q-1}(\M))$
and is called the triple Massey product
$$
\langle [a], [b], [c] \rangle \in
H^{p+q+r-1}(\M)/([a]H^{q+r-1}(\M) + [c]H^{p+q-1}(\M)).
$$
It follows from Proposition 11 that

{\bf Proposition 12.}
{\sl If there is a nontrivial Massey product of classes from 
$H^{\ast}(\M)$, then the algebra $\M$ is not formal.} 

\vskip5mm

\begin{center}
{\bf 5. Cohomology of a blow up of a symplectic manifold}
\end{center}

In this paragraph we expose some computations of cohomologies of
a symplectic blow up, part of which was done in \cite{McDuff}.

Let $(X,\omega)$ be a compact symplectic manifold of dimension $2N$ and
let $Y$ be a symplectic
submanifold of $X$, of dimension $2(N-k)$. Denote by $\widetilde{X}$ 
the blow
up of $X$ along $Y$, by $\pi: \widetilde{X} \rightarrow X$
the projection which is
a degree one mapping, by $V$ the closure of a sufficiently small
tubular neighborhood of $Y$ in $X$ and by
$\widetilde{V} = \pi^{-1}(V)$ the pullback of $V$ by $\pi$.
The boundaries of $\widetilde{V}$ and $V$ are diffeomorphic to a
$S^{2k-1}$-bundle over $Y$ and are diffeomorphic to each other by $\pi$.
By $j$ and $\widetilde{j}$ denote the natural embeddings
$$
j: V \rightarrow X \ \ \mbox{and} \ \ \widetilde{j}: \widetilde{V}
\rightarrow \widetilde{X}.
$$
and by $f$ and $\widetilde{f}$ denote the natural embeddings of pairs
$$
f: (X,\emptyset) \rightarrow (X,V) \ \ \mbox{and} \ \
\widetilde{f}: (\widetilde{X},\emptyset) \rightarrow (\widetilde{X},
\widetilde{V}).
$$

There are the  exact cohomology sequences, for pairs $(\widetilde{X},
\widetilde{V})$
and $(X,V)$, related by the induced homomorphisms $\pi^{\ast}$:
\begin{equation}
\begin{array}{cccccccccc}
\dots & \leftarrow & H^{i+1}(\widetilde{X}) &
\stackrel{\widetilde{f}^{\ast}}{\leftarrow} &
H^{i+1}(\widetilde{X},\widetilde{V}) &
\stackrel{\partial}{\leftarrow} & H^i(\widetilde{V}) &
\stackrel{\widetilde{j}^{\ast}}{\leftarrow} & H^i(\widetilde{X}) &
\stackrel{\widetilde{f}^{\ast}}{\leftarrow} \\
& & \uparrow & & \uparrow & & \uparrow & & \uparrow & \\
\dots & \leftarrow & H^{i+1}(X) &
\stackrel{f^{\ast}}{\leftarrow} & H^{i+1}(X,V) &
\stackrel{\partial}{\leftarrow} & H^i(V) &
\stackrel{j^{\ast}}{\leftarrow} & H^i(X) &
\stackrel{f^{\ast}}{\leftarrow}
\end{array}
\label{exact}
\end{equation}
$$
\begin{array}{cccccccc}
\leftarrow & H^i(\widetilde{X},\widetilde{V}) &
\stackrel{\partial}{\leftarrow} & H^{i-1}(\widetilde{V}) &
\stackrel{\widetilde{j}^{\ast}}{\leftarrow} & H^{i-1}(\widetilde{X}) &
\leftarrow
& \dots . \\
& \uparrow & & \uparrow & & \uparrow & & \\
\leftarrow & H^i(X,V) & \stackrel{\partial}{\leftarrow} & H^{i-1}(V) &
\stackrel{j^{\ast}}{\leftarrow} & H^{i-1}(X) & \leftarrow & \dots .
\end{array}
$$

First, we formulate some simple properties of the diagram of sequences
(\ref{exact}).

{\bf Proposition 13.}
{\sl The vertical homomorphism
$$
\pi^{\ast}: H^i(X,V) \rightarrow H^i(\widetilde{X},\widetilde{V})
$$
is an isomorphism for any $i \geq 0$.}

{\sl Proof.} Denote by $\mbox{Int}\, V$ and $\mbox{Int}\,\widetilde{V}$
the interiors of $V$ and $\widetilde{V}$ which are fibre bundles over $Y$
and notice that by the excision lemma (\cite{Spanier}) 
there are isomorphisms induced by embeddings of pairs
$$
H^i(\widetilde{X},\widetilde{V}) =
H^i(\widetilde{X} \setminus \mbox{Int}\,\widetilde{V},
\widetilde{V} \setminus \mbox{Int}\,\widetilde{V}), \ \
H^i(X,V) =
H^i(X \setminus \mbox{Int}\,V, V \setminus \mbox{Int}\,V).
$$
Notice that 
$\pi^{\ast}: H^{\ast}
(\widetilde{X} \setminus \mbox{Int}\, \widetilde{V}, \widetilde{V}
\setminus
\mbox{Int}\,\widetilde{V}) \rightarrow H^{\ast}(X \setminus \mbox{Int}\, 
V, V \setminus
\mbox{Int}\, V)$ is also an isomorphism and this proves the proposition.

{\bf Proposition 14.}
{\sl For the cohomology groups over the field $F$ of rational, 
real or complex
numbers the vertical homomorphism
$$
\pi^{\ast}: H^i(X;F) \rightarrow H^i(\widetilde{X};F)
$$
is a monomorphism for any $i \geq 0$.}

{\sl Proof.} Without loss of generality we consider the case of real
cohomologies.
Symplectic manifolds are oriented and therefore for any nontrivial
cohomology class $[\tau] \in H^i(X;\R)$ there exists a Poincare dual
class $[\eta] \in H^{2N-i}(X;\R)$ such that
$$
[\tau] \cup [\eta] = [\vol_X] \in H^{2N}(X;\R), \ \
\langle [\vol_X], [X] \rangle = \int_X \vol_X =1.
$$
The homomorphism
$$
\pi^{\ast}: H^{2N}(X;\R) \rightarrow H^{2N}(\widetilde{X};\R)
$$
is a multiplication by $\deg \pi$ and, since $\deg \pi = 1$,
we have
$\pi^{\ast}([\vol_X]) = [\vol_{\widetilde{X}}]$.
Since $\pi^{\ast}$ is a homomorphism of the
cohomology rings, we have
$$
\pi^{\ast}([\tau]) \cup \pi^{\ast}([\eta]) = \pi^{\ast}([\vol_X]) =
[\vol_{\widetilde{X}}]
$$
which implies that $\pi^{\ast}([\tau]) \neq 0$.
This proves the proposition.
  
It is clear from the proof of Proposition 14 that

{\bf Proposition 15.}
{\sl If the cohomology ring $H^{\ast}(X)$ is torsion free then
the vertical homomorphism
$$
\pi^{\ast}: H^i(X) \rightarrow H^i(\widetilde{X})
$$
is a monomorphism for any $i \geq 0$.}

The most important case in which Proposition 15 is applied is $X = 
\C P^N$.

Since $\widetilde{V}$ retracts onto $\widetilde{Y}$ and $V$ retracts onto
$Y$ with preserving the fibration, we derive from Proposition 2
that

{\bf Proposition 16.}
{\sl The vertical homomorphism
$$
\pi^{\ast}: H^i(V) \rightarrow H^i(\widetilde{V})
$$
is a monomorphism for any $i \geq 0$.}

Now we consider consequences
of (\ref{exact}) in some special case.

{\bf Proposition 17.}
{\sl Given a $2(N-k)$-dimensional symplectic submanifold $Y$ of $X = 
\C P^N$,
there are the following short exact sequences:

1) for $i = 2l$ with $0 \leq l \leq (N-k)$,
\begin{equation}
\begin{array}{ccccccccc}
& & & & H^i(\widetilde{V}) & 
\stackrel{\widetilde{j}^{\ast}}{\leftarrow} &
H^i(\widetilde{X}) & \leftarrow & 0 \\
& & &  \partial \swarrow & & & & & \\
0 & \stackrel{\widetilde{f}^{\ast}}{\leftarrow} &
H^{i+1}(\widetilde{X},\widetilde{V}) & & \uparrow \pi^{\ast}
& & \uparrow \pi^{\ast} & & \\
& & & \pi^{\ast} \cdot \partial \nwarrow & & & & & \\
& & & & H^i(V) & \stackrel{j^{\ast}}{\leftarrow} &
\Z = H^i(X) & \leftarrow & 0 \\
\end{array};
\label{short1}
\end{equation}

2) for $i =2l$ with $(N-k+1) \leq l \leq N$,
\begin{equation}
0 \leftarrow H^i(\widetilde{V}) 
\stackrel{\widetilde{j}^{\ast}}{\leftarrow}
H^i(\widetilde{X}) \stackrel{\pi^{\ast}f^{\ast}}{\leftarrow} 
\Z = H^i(X) = H^i(X,V) \stackrel{\partial}{\leftarrow} 0;
\label{short2}
\end{equation}

3) for $i = 2l+1$ with $i+1 \leq \dim Y = 2(N-k)$,
\begin{equation}
\begin{array}{ccccccccc}
0 & \leftarrow & H^{i+1}(\widetilde{X},\widetilde{V}) &
\stackrel{\partial}{\leftarrow} & H^i(\widetilde{V}) &
\stackrel{\widetilde{j}^{\ast}}{\leftarrow} & H^i(\widetilde{X})
& \leftarrow & 0 \\
& & & \pi^{\ast}\partial \nwarrow & \uparrow \pi^{\ast}& & & & \\
& & & & H^i(V) & \leftarrow & 0 & & \\
\end{array};
\label{short3}
\end{equation}

4) for $i = 2l+1$ with $i+1 > \dim Y = 2(N-k)$,
\begin{equation}
0 \leftarrow \pi^{\ast}(H^{i+1}(X)) = \Z
\leftarrow  H^{i+1}(\widetilde{X},\widetilde{V})
\stackrel{\partial}{\leftarrow}  H^i(\widetilde{V}) 
\stackrel{\widetilde{j}^{\ast}}{\leftarrow}  H^i(\widetilde{X})
\leftarrow 0.
\label{short4}
\end{equation}
}

{\sl Proof.}

1) Since the embedding $Y \subset V \subset X$ is symplectic,
$j^{\ast}:H^i(X) \rightarrow H^i(V) = H^i(Y)$ is a monomorphism for
any $i \leq 2(N-k) + 1$. This implies that
\begin{equation}
f^{\ast} = 0 \ \ \mbox{for \ \ $i \leq 2(N-k) + 1$}.
\label{fact1}
\end{equation}
Since the diagram (\ref{exact}) is commutative and
$\pi^{\ast}:H^i(X,V) \rightarrow H^i(\widetilde{X},\widetilde{V})$
is an isomorphism, we have
\begin{equation}
\tilde{f}^{\ast} = 0 \ \ \mbox{for \ \ $i \leq 2(N-k) + 1$}.
\label{fact2}
\end{equation}
Now (\ref{exact}) together with (\ref{fact1}) and (\ref{fact2}) implies
(\ref{short1}).

2) For $ i =2l > 2(N-k)$,
$H^i(V) = H^{i-1}(V) = 0$ and $H^{i+1}(X) = 0$.
Therefore the lower exact sequence in (\ref{exact}) together with
Proposition 13 implies
\begin{equation}
H^{i+1}(X,V) = H^{i+1}(\widetilde{X},\widetilde{V}) = 0 \ \
\mbox{for $i =2l > 2(N-k)$}
\label{fact3}
\end{equation}
and
\begin{equation}
H^i(X) \stackrel{f^{\ast}}{\approx} H^i(X,V) \ \
\mbox{for $i =2l > 2(N-k)$}.
\label{fact4}
\end{equation}
Now (\ref{exact}) together with (\ref{fact3}) and (\ref{fact4}) implies
(\ref{short2}).

3) For $i = 2l+1$, $H^i(X) = 0$.
Together with (\ref{fact2}) for $i +1 = 2l+2 \leq \dim Y$
this implies (\ref{short3}).

4) For $i = 2l+1$    with $i+1 > \dim Y$, we have
$H^{i+1}(V) = 0$ and $H^i(X) = 0$. The first equality implies that
the homomorphism $\widetilde{j}^{\ast}:H^{i+1}(\widetilde{X}) \rightarrow
H^{i+1}(\widetilde{V})$ is trivial: $\widetilde{j}^{\ast} = 0$. 
The second equality together with Proposition 13 and 
the commutativity of (\ref{exact})
implies that the homomorphism 
$\widetilde{f}^{\ast}:H^i(\widetilde{X},\widetilde{V}) \rightarrow  
H^i(\widetilde{X})$ is also trivial. In this event 
the corresponding fragment
of (\ref{exact}) reduces to (\ref{short4}).

Proposition 17 is proved.

Recall that the minimal model of $\C P^m$ is a differential graded
algebra $\M_{\C P^m}$ freely generated by elements $x$ and $y$,
with $\deg x =2$ and $\deg y = 2k-1$, and 
the differential acts as follows:
$dx =0, dy = x^k$. 

The minimal model of $\widetilde{Y}$ is computed by using
Proposition 2 as follows.

{\bf Proposition 18.}
{\sl Given  the minimal model  $\M_Y$ for $Y$ and a 
$U(k)$-vector bundle
$E \rightarrow Y$, the minimal model $\M_{\widetilde{Y}}$
for the projectivization $\widetilde{Y}$ is isomorphic to
$$
\M_{\widetilde{Y}} = \M_Y \otimes_d \M_{\C P^{k-1}}
$$
which is a differential graded algebra freely generated by elements of
$\M_Y$ and $\M_{\C P^{k-1}}$ and the differential $d$ acts on
$\M_{\widetilde{Y}}$ as follows

1) its restriction onto $\M_Y$ equals the differential of
$\M_Y$;

2) $dx = 0$ and
$$
dy = x^k + c_1 x^{k-1} + \dots + c_{k-1} x + c_k
$$
where the elements $c_j \in \M_Y$ represent the rational
Chern classes $c_j(E)$ via the isomorphism
$H^{\ast}(\M_Y) = H^{\ast}(Y;\Q)$.}

This proposition is quite clear and also follows
from general facts on relations between minimal models 
and Serre fibrations (\cite{Thomas}).

\vskip5mm

\begin{center}
{\bf 6. Examples of nonformal simply connected symplectic manifolds}
\end{center}

{\bf A) A series of nonformal symplectic nilmanifolds.}

Consider the algebra $W(1)$ of formal vector fields on the line.
This is a topological infinite-dimensional algebra for which a basis
is given by linear differential operators
$$
e_k = x^{k+1} \frac{d}{dx}, \ k=-1,0,1,\dots,
$$
the Lie brackets in this basis takes the form
\begin{equation}
[e_i,e_j] = (j-i)e_{i+j}, \ \ i,j \geq -1.
\label{liebr1}
\end{equation}
The algebra $W(1)$ has a natural filtration
$$
\dots \subset {\cal L}_1(1) \subset {\cal L}_0(1) \subset
{\cal L}_{-1}(1) \subset W(1)
$$
where ${\cal L}_k(1)$ is the subalgebra spanned by $e_k,e_{k+1},\dots$.
Here we use the notations from \cite{Fuks} which are widely accepted.
We also refer to this book for information about the cohomologies of 
the Lie algebras
${\cal L}_k(1)$ and their relations to other subjects.

In the sequel we consider a series of finite-dimensional nilpotent Lie
algebras
$$
{\cal V}_n = {\cal L}_1(1)/{\cal L}_{n+1}(1), \ n=3,4,\dots.
$$
We denote by $V_n$ the corresponding Lie groups. Hence ${\cal V}_n$
is an $n$-dimensional algebra with  basis $\{e_1,\dots,e_n\}$
and the Lie brackets
\begin{equation}
[e_i,e_j] =
\cases{
(j-i)e_{i+j}, & for $i+j \leq n$ \cr
0, & for $i+j>n$.
}
\label{liebr2}
\end{equation}
The structure constants of ${\cal V}_n$ are rational and therefore
${\cal V}_n$ possesses uniform lattices. We take one of them which may be
called canonical (\cite{Maltsev}). The group $V_n$ is isomorphic to
$({\cal V}_n,\times)$ where the multiplication $\times$ is given
by the Campbell--Hausdorff formula.
The basis $\{e_1,\dots,e_n\}$
multiplicatively generates a subgroup $\Gamma_n$ and we obtain 
as a result an infinite family of finite-dimensional nilmanifolds
$$
M(n) = V_n/\Gamma_n, \ n=3,4,\dots.
$$

The group $V_3$ is the Heisenberg group ${\cal H}$,
$M(3)$ is ${\cal H}/{\cal H}_Z$, and $M(3) \times S^1$ is
the Kodaira--Thurston manifold $\widetilde{M}$.

There are only three $4$-dimensional nilpotent groups: 1) the commutative
group $\R^4$ which admits a left-invariant K\"ahler structure;
2) $V_3 \oplus \R$ which admits left-invariant complex and symplectic
structures but no left-invariant K\"ahler structure;
3) $V_4$, a three step nilpotent group, which admits a left-invariant
symplectic structure but no left-invariant complex structure (\cite{FGG}).

Let $\{\omega_1,\dots,\omega_n\}$ be a basis for left-invariant $1$-forms,
on $V_n$, dual to the basis $\{e_1,\dots,e_k\}$. Then (\ref{liebr2}) implies
(see also Section 4)
\begin{equation}
d \omega_k = (k-2) \omega_1 \wedge \omega_{k-1} + (k-4) \omega_2 \wedge
\omega_{k-2} + \dots.
\label{diff2}
\end{equation}
By Proposition 9, we have

{\bf Proposition 19.}
{\sl The minimal model for $M(n) = V_n/\Gamma_n$ is as follows:}
$$
\M(n) = (\Lambda(x_1,\dots,x_n),d) \ \mbox{with} \ \deg x_k =1 \
\mbox{for} \ k=1,\dots,n,
$$
$$
dx_1= dx_2=0,\
dx_k = (k-2) x_1 \wedge x_{k-1} + (k-4) x_2 \wedge
x_{k-2} + \dots \ \mbox{for} \ k \geq 3.
$$

In fact, the algebra $\M(n)$ is bigraded with the second grading given by
$\deg' x_i = i$. In this event the bidegree of the
differential is $(1,0)$.

{\bf Proposition 20.}
{\sl The form
$$
\Omega_{2m} = (2m-1) \omega_1 \wedge \omega_{2m} + (2m-3) \omega_2 \wedge
\omega_{2m-1} + \dots + \omega_m \wedge \omega_{m+1}
$$
is a left-invariant symplectic form on $V_{2m}$ for $m \geq 2$.}

The proof of this proposition is as follows. Take an extension of
${\cal V}_{2m}$ by adding a new generator $\omega_{2m+1}$ such that
$d\omega_{2m+1} = \Omega_{2m}$. It follows from (\ref{diff2}) that
this would be ${\cal V}_{2m+1}$ and hence 
$d^2 \omega_{2m+1} = d\Omega_{2m} =
0$. This proves the proposition.

Since $\Omega_{2m}$ is invariant, its pushforward on the quotient space
$V_{2m}$ is well-defined and we preserve for the pushfoward 
the same notation.
It is clear that  the integrals of $\Omega_{2m}$ over cycles in
$H_2(V_{2m};\Z)$ are integer.
This implies

{\bf Corollary 1.}
{\sl The nilmanifolds $M(2m)$ admit integer symplectic forms.}

We will not compute the cohomologies of these nilmanifolds and
restrict ourselves to the following fact.

{\bf Proposition 21.}
{\sl For any $m \geq 2$,

1) $H^1(M(2m);\Q)$ is spanned by $[x_1]$ and $[x_2]$;

2) $[x_2 \wedge x_3] \neq 0$ in $H^2(M(2m);\Q)$.}

Since $\M(2m)$ contains only two closed generators, the first statement is
clear. Notice now that $x_2 \wedge x_3$ has bidegree $(2,5)$ and recall that
$d$ has bidegree $(1,0)$. If $x_2 \wedge x_3 = du$, then $u$ has to be
proportional to $x_5$ but, by (\ref{diff2}), $dx_5 = 3x_1 \wedge x_4 +
x_2 \wedge x_3$. The proposition is proved.

{\bf Theorem 1.}
{\sl The symplectic nilmanifolds $M(2m)$ are nonformal.}

Taking into account that $dx_3 = x_1 \wedge x_2$ in $\M(2m)$,
a proof of the theorem  is obtained from Proposition 21 by the
same reasonings as the proof of Proposition 10 above.

{\bf B) Nonformal simply connected symplectic manifolds.}

By Proposition 1, there exist symplectic embeddings of
$\widetilde{M} = M(3) \times S^1$ into $\C P^N$ with $N \geq 5$ and
of $M(2m)$ into $\C P^N$ with $N \geq 2m+1$ such that
the symplectic form (\ref{thurston}) and $\Omega_{2m}$ are the pullbacks of
the Fubiny--Study forms $\omega_{FS}$ on $\C P^N$ under these
embeddings.

Now, realizing these nilmanifolds as symplectic submanifolds of
complex projective spaces,
denote by $\widetilde{\C P}^N$ the symplectic blow up of $\C P^N$
along $\widetilde{M}$ and denote by
$\widetilde{X}_m(N)$ the symplectic blow up of $\C P^N$
along $M(2m)$.

{\bf Theorem 2.}
{\sl For $m \geq 2$ and $N \geq 2m+1$ the symplectic manifolds
$\widetilde{\C P}^N$ and $\widetilde{X}_m(N)$ are simply connected and
nonformal.}

{\sl Proof of Theorem 2.}

By Proposition 3, these manifolds are simply connected.

We will use the same notation as in Section 5 and denote by $Z$ 
the closure of the complement to
$\widetilde{V}$ and $V$ in $X$ and $Y$. We use the following notation for
the embeddings
$$
Z \stackrel{\widetilde{i}_1}{\longleftarrow} \partial V
\stackrel{\widetilde{i}}{\longrightarrow} \widetilde{V},
\ \
Z \stackrel{\widetilde{j}_1}{\longrightarrow} \widetilde{X}
\stackrel{\widetilde{j}}{\longleftarrow} \widetilde{V}
$$
and recall that $\widetilde{j}$ was already introduced in Section 5 and
$\widetilde{i}$ was introduced in Section 3.
We also consider rational cohomologies.

Take the exact Mayer--Vietoris sequence for the pair $(Z,\widetilde{V})$:
$$
\dots \rightarrow H^q(Z \cup \widetilde{V})
\stackrel{\widetilde{j}_1^{\ast} \oplus
\widetilde{j}^{\ast}}{\longrightarrow}
H^q(Z) \oplus H^q(\widetilde{V})
\stackrel{\widetilde{i}_1^{\ast} - \widetilde{i}^{\ast}}{\longrightarrow}
 H^q(Z \cap \widetilde{V}) \rightarrow \dots.
$$
By (\ref{short1}),(\ref{short3}), and (\ref{short4}),
$\widetilde{j}^{\ast}$ is a monomorphism for
$q \leq \dim Y +1 $ and we derive that there are the following splittings:
\begin{equation}
0 \rightarrow H^q(\widetilde{X})
\stackrel{\widetilde{j}_1^{\ast} \oplus
\widetilde{j}^{\ast}}{\longrightarrow}
H^q(Z) \oplus H^q(\widetilde{V})
\stackrel{\widetilde{i}_1^{\ast} - \widetilde{i}^{\ast}}{\longrightarrow}
H^q(\partial \widetilde{V}) \rightarrow 0 \ \ \mbox{for $q \leq \dim Y +1$}.
\label{vietoris}
\end{equation}
 
For any $Y$ of the form $M(2m)$ or $M(3)\times S^1$ its minimal model
contains closed generators $x_1$ and $x_2$ and a generator $x_3$ such that
$d x_3 = x_1 \wedge x_2$ (see (\ref{thurston}) and Propositions 9 
and 19).

Take elements $(0,a \cup [x_1])$ and $(0,a \cup [x_2])$ in
$H^3(Z) \oplus H^3(\widetilde{V})$ and notice that by Proposition 2,
$\widetilde{i}^{\ast}(a) = 0$ which together with (\ref{vietoris})
implies that there are $u_1, u_2 \in H^3(\widetilde{X})$ such that
$\widetilde{j}^{\ast}(u_k) = a \cup [x_k]$ for $k =1,2$.

We are left to prove two lemmas.

{\bf Lemma 1.}
{\sl For $m \geq 3$ the triple Massey product $\langle u_2,u_1,u_2 \rangle$
is defined and nontrivial in
$H^8(\widetilde{X}_m(N))/u_2 \cup  H^5(\widetilde{X}_m(N))$.}

{\bf Lemma 2.}
{\sl For the symplectic manifolds $\widetilde{\C P}^N$ and
$\widetilde{X}_2(N)$
the triple Massey product $\langle u_2, v, u_2 \rangle$ is defined
where $v = \pi^{\ast}([\omega])$ and $\omega$ is a symplectic form on
$X = \C P^N$. These products are also nontrivial in
$H^7/(u_2 \cup H^4)$.}

{\sl Proof of Lemma 1.}

Since $[x_1 \wedge x_2] = 0$ in $M_{\widetilde{V}}$,
we have
$$
\widetilde{j}^{\ast}(u_1 \cup u_2) =\widetilde{j}^{\ast}(u_1) \cup
\widetilde{j}^{\ast}(u_2) = a^2 \cup ([x_1] \cup [x_2]) = 0.
$$
But $\deg (u_1 \cup u_2) = 6 \leq \dim M(2m)$ and, since, by
(\ref{short1}), $\widetilde{j}^{\ast}:H^6(\widetilde{X}) \rightarrow
H^6(\widetilde{V})$ is a monomorphism, $u_1 \cup u_2 = 0$ in
$H^6(\widetilde{X})$. Therefore the triple product
$\langle u_2, u_1, u_2 \rangle$ is defined.

The image of the triple product in $H^8(\widetilde{V})$ equals
$ a^3 \cup ([x_3] \cup [x_2]) \ \mbox{modulo} \ (a \cup [x_2]) \cup
H^5(\widetilde{V})$. It is easily computed by using Proposition 18 that
$$
a^3 \cup ([x_3] \cup [x_2]) \neq 0 \ \ \mbox{modulo} \ (a \cup [x_2]) \cup
H^5(\widetilde{V}).
$$
We have $\widetilde{j}^{\ast}(u_2 \cup H^5(\widetilde{X})) \subset
(a \cup [x_2]) \cup H^5(\widetilde{V})$ and 
this implies that the triple product
$\langle u_2, u_1, u_2 \rangle$  is also nontrivial modulo
$u_2 \cup H^5(\widetilde{X})$ and proves the lemma.

{\sl Proof of Lemma 2.}

By Proposition 18, $a = [x] \in H^2(\widetilde{V})$.

Compute
$$
\widetilde{j}(u_2 \cup v) = [x \wedge x_2] \cup [ A x_1 \wedge x_4 +
x_2 \wedge x_3] =
$$
$$
A [x \wedge x_2 \wedge x_1 \wedge x_4] =
-A [d(x \wedge x_3 \wedge x_4)]
$$
where $A = 1$ for $Y = \widetilde{M}$ and $A = 3$ for
$Y = M(4)$.
Therefore $\widetilde{j}^{\ast}(u_2 \wedge v) = 0$.

By (\ref{short3}), $\widetilde{j}^{\ast}:H^5(\widetilde{X}) \rightarrow
H^5(\widetilde{V})$ is a monomorphism and therefore
$$
u_2 \cup v = 0.
$$
Hence the triple product $\langle u_2, v, u_2 \rangle$ is defined
and its image in $H^7(\widetilde{V})$ equals
$-A a^2 \cup [x_2 \wedge x_3 \wedge x_4] \ \mbox{modulo} \
(a \cup [x_2]) \cup
H^4(\widetilde{V})$. It is easily computed by using Proposition 18 that
$$
a^2 \cup [x_2 \wedge x_3 \wedge x_4] \neq 0 \ \ \mbox{modulo} \
(a \cup [x_2]) \cap H^5(\widetilde{V}).
$$
Since $\widetilde{j}^{\ast}(u_2 \cup H^4(\widetilde{X})) \subset
(a \cup [x_2]) \cup H^4(\widetilde{V})$,
this implies that the triple product
$\langle u_2, v, u_2 \rangle$  is nontrivial modulo
$u_2 \cup H^4(\widetilde{X})$ and proves the lemma.

Now Theorem 2 follows from Lemmas 1 and 2 and Proposition 12.

{\bf Proposition 22.}
{\sl Let $X_1$ and $X_2$ be simply connected manifolds of dimension $N$
such that there is a nontrivial triple Massey product in $H^q(X_1)$
with $q \leq N-3$. Then there is a nontrivial Massey product in
$H^q(X_1 \# X_2)$.}

 Note that one can prove a general proposition: {\sl The connected sum of two
simply connected manifolds is formal if and only if each summand is formal},
but actually we need the slightly weaker version stated above.

{\sl Proof.} Take the Mayer--Vietoris sequence for the pair
$(\overline{X_1 \setminus D},\overline{X_2 \setminus D})$ where $D$ is an
$N$-disc and the overline stands for closure.
Since $H^q(S^{N-1}) = H^q(\overline{X_1 \setminus D} \cap
\overline{X_2 \setminus D}) = 0$ for $1 \leq q \leq N-2$, the embeddings
$j_1: \overline{X_1 \setminus D} \rightarrow X_1 \# X_2$ and
$j_2: \overline{X_2 \setminus D} \rightarrow X_1 \# X_2$ induce
isomorphisms
$$
H^q(X_1 \# X_2)
\stackrel{j_1^{\ast} \oplus j_2^{\ast}}{\longrightarrow}
H^q(\overline{X_1 \setminus D}) \oplus
H^q(\overline{X_2 \setminus D})
$$
for $1 \leq q \leq N-3$.

Considering the Mayer--Vietoris sequence for the pair
$((\overline{X_1 \setminus D},\overline{D})$ we conclude that there are
isomorphisms
$$
H^q(X_1) \rightarrow H^q(\overline{X_1 \setminus D})
$$
for $1 \leq q \leq N-3$.

Therefore, if there is a nontrivial triple Massey product of degree 
$\leq N-3$
in the cohomologies of $X_1$ it survives in the cohomologies of
$X_1 \# X_2$. This proves the proposition.

{\bf Corollary 2.}
{\sl For any $k \geq 1$ the symplectic manifolds
$\widetilde{\C P}^N \#\,k\overline{\C P}^N$ and
$\widetilde{X}_m(N)\#\,k\overline{\C P}^N$ are nonformal.}

The proof is clear, because the nontrivial triple products in
$\widetilde{\C P}^N$ and $\widetilde{X}_2(N)$, constructed in the proof of
Lemma 2, are
of degree $7$ and the nontrivial triple products in
$\widetilde{X}_m(N)$, with $m \geq 3$, constructed in the proof  
of Lemma 1, are of degree $8$ and therefore all of them survives
symplectic blows up at points.

\vskip5mm

\begin{center}
{\bf 7. Some remarks on formality and symplectic Hodge theory}
\end{center}

First recall some definitions from symplectic Hodge theory.

Take, for simplicity, local coordinates $\{x^1,\dots,x^{2N}\}$ on a
manifold with a symplectic form
$$
\omega = \omega_{ij} dx^i \wedge dx^j.
$$
and take the inverse tensor to $\omega_{ij}$, defined by
$$
\omega^{ij}\omega_{jk} = \cases{1 & if $i = k$ \cr 0 & otherwise}.
$$
The symplectic inner product of differential $k$-forms is
$$
\langle \alpha, \beta \rangle = \omega^{i_1 j_1} \dots \omega^{i_k j_k}
\alpha_{i_1 \dots i_k} \beta_{j_1 \dots j_k}.
$$
Now define the symplectic star operator
by the condition
$$
\alpha \wedge \ast \beta  = \frac{1}{n!} \langle \alpha, \beta \rangle
\, \omega^n.
$$
The latter definition enables us to introduce a codifferential $\delta$
of degree $-1$:
$$
\delta = \ast \cdot d \cdot \ast, \ \ \ \delta^2 = 0.
$$
Finally call a differential form $\alpha$ symplectically harmonic if
$$
d \alpha = \delta \alpha = 0.
$$
These definitions are due to Koszul (\cite{Koszul}) and Brylinski
(\cite{Brylinski}).

Brylinski proved that on a compact K\"ahler manifold $M$ 
each cohomology class
in $H^{\ast}(M;\C)$ is realized by a symplectically harmonic form
and he also conjectured that this analog of the Hodge theorem for Riemannian
manifolds is valid in general (\cite{Brylinski}).

Mathieu proved that a compact symplectic manifold satisfies the Brylinski
conjecture if and only if it satisfies the Hard Lefschetz condition
(\cite{Mathieu}). This theorem implies that the Kodaira--Thurston manifold
$\widetilde{M}$ and the McDuff manifold $\widetilde{\C P}^5$
serve as counterexamples to the Brylinski conjecture.

Recently Merkulov (\cite{Merkulov})  considered the problem of formality
for simpy connected symplectic manifolds satisfying the Hard Lefschetz
condition. His conception of formality deals with quasi-isomorphisms
which are additive only and do not respect the multiplicative structures
of cochain complexes, so it has no a real topological meaning.   

Nevertheless we introduce the following

{\bf Conjecture 1.}
{\sl A simply connected compact symplectic manifold $M$ is formal if and
only if it satisfies the Hard Lefschetz condition, or, which is equivalent,
each cohomology class in $H^{\ast}(M)$ is realized by a symplectically
harmonic form.}

We also would like to state

{\bf Conjecture 2.}
{\sl Let $Y \rightarrow X$ be a symplectic embedding of a nonformal
simply connected manifold into a simply connected compact manifold $X$.
Then the blow up of $X$ along $Y$ is nonformal.}

\vskip1cm

{\bf Final remarks.}

This paper was started during the visit of the second author (I.A.T.)
to the University of Montpellier and was finished during his visit to
SFB 288 in the Technical University of Berlin. The authors were 
also partially
supported by the Russian Foundation for Basic Researches (grants
96-01-00182a (I.K.B.) and 96-15-96877 and 98-01-00749 (I.A.T.)).


\vskip15mm


\begin{thebibliography}{MMM}


\bibitem{BT}
Babenko, I. K., Taimanov, I. A.,
On existence of nonformal simply connected symplectic manifolds,
Uspekhi Mat. Nauk {\bf 53}:5 (1998), 225--226 (Russian);
English translation: Russian Math. Surveys {\bf 53}:5 (1998).


\bibitem{BG}
Benson, C., Gordon, C.,
K\"ahler and symplectic structures on nilmanifolds,
Topology {\bf 27} (1988), 513--518.


\bibitem{Brylinski}
Brylinski, J. L.,
A differential complex for Poisson manifolds,
J. Differ. Geom. {\bf 28} (1988), 93--114.


\bibitem{CFG}
Cordero, L. A., Fernandez, M., Gray A.,
Symplectic manifolds with no K\"ahler structure,
Topology {\bf 25} (1986), 375--380.


\bibitem{DGMS}
Deligne, P., Griffiths, P., Morgan, J., Sullivan, D.,
Real homotopy theory of K\"ahler manifolds,
Invent. Math.  {\bf 19} (1975), 245--274.



\bibitem{FGG}
Fernandez, M., Gotay, M., Gray, A.,
Compact parallelizable four dimensional symplectic and complex manifolds,
Proc. Amer. Math. Soc. {\bf 103} (1988), 1209--1212.


\bibitem{Fuks}
Fuks, D. B.,
Cohomology of infinite-dimensional Lie algebras,
Consultants Bureau, New York, 1986.


\bibitem{Gompf}
Gompf, R. E.,
A new construction of symplectic manifolds,
Ann. of Math. (II) {\bf 142} (1995), 527--595.


\bibitem{GM}
Griffiths, P., Morgan, J.,
Rational Homotopy Theory and Differential Forms,
Birkh\"auser, Basel, 1981.


\bibitem{Nice}
Gromov, M. L.,
A topological technique for the construction of solutions of differential
equations and inequalities,
Actes Congr\'es Intern. Math. (Nice, 1970), Gauthier-Villars, Paris,
No. 2, 1971, 221--225.


\bibitem{Gromov}
Gromov, M. L.,
Partial Differential Relations, Springer, Berlin-Heidelberg, 1986.


\bibitem{Husemoller}
Husemoller, D.,
Fibre bundles, McGraw-Hill, New York, 1966.


\bibitem{Koszul}
Koszul, J. L.,
Crochet de Schouten--Nijenhuis et cohomologie,
in ``Elie Cartan et les math\'ematiqeu d'aujourd'huis'',
Asterisque, 1985, 251--271.


\bibitem{LO}
Lupton, G., Oprea, J.,
Symplectic manifolds and formality,
J. Pure Appl. Algebra {\bf 91} (1994), 193--207.

\bibitem{Maltsev}
Mal'tsev, A. I.,
On a class of homogeneous spaces,
Izvestia Akad. Nauk SSSR Ser. Mat.
{\bf 3} (1949), 9--32 (Russian); English translation:
Amer. Math. Soc. Transl. (1) {\bf 9} (1962), 276--307.


\bibitem{Mathieu}
Mathieu, O.,
Harmonic cohomology classes of symplectic manifolds,
Comment. Math. Helv. {\bf 70} (1995), 1--9.


\bibitem{McDuff}
McDuff, D.,
Examples of symplectic simply connected manifolds 
with no K\"ahler structure,
J. Differ. Geom. {\bf 20} (1984), 267--277.


\bibitem{MS}
McDuff, D., Salamon, D.,
Introduction to Symplectic Topology,
Clarendon Press, Oxford, 1995.


\bibitem{Merkulov}
Merkulov, S. A.,
Formality of canonical symplectic complexes and Frobenius manifolds,
Internat. Math. Res. Notices {\bf 14} (1998), 727--733.

\bibitem{NM}
Neisendorfer, J., Miller, T.,
Formal and coformal spaces,
Illinois J. Math. {\bf 22} (1978), 565--579.


\bibitem{Nomizu}
Nomizu, K.,
On the Cohomology of Homogeneous Spaces of Nilpotent Lie Groups,
Ann. of Math. (II) {\bf 59} (1954), 531--538.



\bibitem{Spanier}
Spanier, E. H.,
Algebraic Topology,
McGraw-Hill, New York, 1966.


\bibitem{Sullivan}
Sullivan, J.,
Infinitesimal computations in topology,
Publ. IHES {\bf 47} (1978), 269--331.


\bibitem{Thomas}
Thomas, J. C.,
Rational homotopy of Serre fibrations,
Annales Inst. Fourier {\bf 31}:3 (1981), 71--90.


\bibitem{Thurston}
Thurston, W.,
Some simple examples of compact symplectic manifolds,
Proc. Amer. Math. Soc. {\bf 55} (1976), 467--468.


\bibitem{Tischler}
Tischler, D.,
Closed $2$-forms and an embedding theorem for symplectic manifolds,
J. Differ. Geom. {\bf 12} (1977), 229--235.


\bibitem{TO}
Tralle, A., Oprea, J.,
Symplectic Manifolds with no K\"ahler Structure,
Lect. Notes in Math. {\bf 1661}, Springer, Berlin Heidelberg, 1997.


\end{thebibliography}
\end{document}